\documentclass[12pt,oneside,CJK,reqno]{amsart}
\usepackage{mathrsfs}
\usepackage{bbm}
\usepackage{amssymb}
\usepackage{enumerate,amsmath,amssymb,amsthm}
\usepackage{CJK}
\usepackage{amsfonts}
\newcommand{\be}{\begin{eqnarray}}
\newcommand{\ee}{\end{eqnarray}}
\newcommand{\ce}{\begin{eqnarray*}}
\newcommand{\de}{\end{eqnarray*}}
\newtheorem{theorem}{Theorem}[section]
\newtheorem{lemma}[theorem]{Lemma}
\newtheorem{remark}[theorem]{Remark}
\newtheorem{definition}[theorem]{Definition}
\newtheorem{proposition}[theorem]{Proposition}

\newtheorem{corollary}[theorem]{Corollary}
\def\e{\varepsilon}
\def\t{\theta}
\def\a{\alpha}
\def\om{\omega}
\def\b{\beta}
\def\p{\partial}

\def\d{\delta}

\def\[{{\Big[}}
\def\]{{\Big]}}
\def\<{{\langle}}
\def\>{{\rangle}}
\def\({{\Big(}}
\def\){{\Big)}}

\def\no{\nonumber}
\def\bt{\begin{theorem}}
\def\et{\end{theorem}}
\def\bl{\begin{lemma}}
\def\el{\end{lemma}}
\def\br{\begin{remark}}
\def\er{\end{remark}}
\def\bd{\begin{definition}}
\def\ed{\end{definition}}
\def\bp{\begin{proposition}}
\def\ep{\end{proposition}}
\def\bc{\begin{corollary}}
\def\ec{\end{corollary}}

\def\cF{{\mathcal F}}

\def\cT{{\mathcal T}}

\def\sB{{\mathscr B}}

\def\sF{{\mathscr F}}

\def\sS{{\mathscr S}}

\def\sV{{\mathscr V}}
\def\sW{{\mathscr W}}

\def\mE{{\mathbb E}}

\def\mN{{\mathbb N}}

\def\mR{{\mathbb R}}

\def\mW{{\mathbb W}}

\def\geq{\geqslant}
\def\leq{\leqslant}

\def\Om{\Omega}

\def\dotw{\dot{w}_n(t)}
\def\dotwu{\dot{w}_n(u)}

\def\s{\sigma}
\def\vph{\varphi}
\def\pf{\begin{proof}}
\def\ef{\end{proof}}

\def\pvph{\partial\vph}

\allowdisplaybreaks

\begin{document}
\begin{CJK}{GBK}{song}

\title{Support theorem for stochastic variational inequalities*}
\author{Jiagang Ren, Siyan Xu}

\subjclass{}

\date{}
\dedicatory{ School of Mathematics and Computational Science, Sun Yat-Sen University,\\
 Guangzhou, Guangdong, 510275, P.R.China\\
Emails: J. Ren: renjg@mail.sysu.edu.cn, S. Xu: xsy\_00@hotmail.com}

\thanks{*Work supported by NSFC (Grant no. 10871215), Ph.D. Programs Foundation of Ministry
of Education of China (Grant no. 20060558051)
and China  postdoctoral science foundation (Grant no. 20090460824).}

\begin{abstract}
We prove a support theorem of the type of
Stroock-Varadhan for solutions of stochastic variational inequalities.\\\\
{\bf Keywords:} {\it stochastic variational inequality,
multivalued maximal monotone operator, limit theorems,
approximate continuity, support theorem.}\\\\
{\bf AMS Mathematics Subject Classification (2000):} 60H10, 60F99, 47H04. \\
\end{abstract}
\maketitle

\section{Introduction and main result}
The aim of this paper is to describe
the support of the law of the solution of the following
stochastic variational inequality (SVI in short):
\be\label{eq0}
\left\{
\begin{array}{ll}
dX(t)\in b(X(t))dt+\sigma(X(t))\circ dw(t)-\pvph(X(t))dt, ~~t\in \mR_+,\\
X(0)=x\in \overline{D(\pvph)},
\end{array}
\right.
\ee
where $\vph$ is a convex function and $\pvph$ is its subdifferential. SVIs of this type have been investigated by many authors in the past two decades (see e.g.
\cite{ba1, ben, ce, ce2, ce3} and reference therein) and they include as a special case
stochastic differential equations (SDE in short) in convex domains reflected at the
boundaries. Also it is needless to say that they reduce to usual SDEs if $\vph$ is differentiable.
For connection of SVIs with parabolic and elliptic Neumann problems and parabolic variational inequalities we refer to \cite{ba1, par}.

The support of the law of a diffusion defined by an SDE was first characterized
by Stroock and Varadhan in \cite{st} and this result has undergone various
extensions ever since, among which for the most recent ones we
mention only \cite{nak}.  For reflected diffusions in {\it smooth}
($C^2$) domains this was done by Doss and Priouret in \cite{do}.
Their approach, however, involves a heavy localization procedure
which does not seem applicable to convex domains with only Lipschitz
boundaries.

As was originally done by Stroock and Varadhan and is now a standard approach,
it consists of proving two inclusion relations to
characterize the support of a diffusion: the direct one
and the inverse one. The direct inclusion involves essentially a limit theorem
for the equation. Such theorem was proved in
\cite{do} for reflected diffusions in  {\it smooth domains} and
in \cite{pet} for reflected diffusions
in convex domains but with {\it constant diffusion coefficients} and some other extra
assumptions which are not easy to verify. In Section 3, we shall prove a limit theorem
for general SVIs in the form of (1). Here we would like to point out that compared
with our recent work \cite{rx}, the main difference is that in that paper we approximate
SVIs by ODEs while in the present one we approximate SVIs by ordinary variational
inequalities. Each one has its advantage and its disadvantage: the former is good for
establishing various regularity properties for solutions of SVIs and the latter is
adequate for determining the support.

The inverse inclusion is deduced from the Denjoy approximate
continuity of solutions of SVIs at sufficiently regular sample
paths. For SDEs with {\it smooth} reflecting boundaries this was
proved in \cite{do} and for multivalued SDEs  with {\it bounded}
multivalued maximal monotone operators this was proved recently in
\cite{xu}. This boundedness assumption is, however, so strong that it even excludes the case of
reflected diffusions. In the present paper we shall be able to remove this
assumption, see Theorem \ref{thsection4}.

Combining the two inclusions will yield the main result of the paper,
Theorem \ref{main}.

The paper is organized as follows: in Section 2 we prepare necessary preliminary materials
and we state our main result in Section 3.
Section 4 and Section 5 will be devoted to the proof of the main result.
\section{Preliminaries}
In this section we collect some materials which will be needed below.

$\cT$ will be the space of c\`adl\`ag real valued functions of
finite variation defined on $\mR_+$ with the metric
$$
d(f,g):=\int_0^\infty \frac{|f(t)-g(t)|}{1+|f(t)-g(t)|}e^{-t}dt.
$$
This metric obviously corresponds to convergence in Lebesgue measure and is also equivalent to
the convergence at each point of continuity of the target function.
It is trivial that $(\cT,d)$ is separable and it is an easy consequence
of the standard diagonalization argument (see, e.g., \cite[p.210-212]{ca})
that it is also complete. Hence it is a Polish space.
Furthermore, we have (see e.g. \cite[Lemma 13.15]{ca})
\bl\label{com}
For any increasing positive function $C(t)$, the set
$$
\{f: |f|_t\leq C(t), \forall  t\}
$$
is compact in $(\cT,d)$, where $|f|_t$ stands for the total variation of $f$ on $[0,t]$.
\el

For $\kappa\in \cT$ which is increasing, define its inverse by
$$
\kappa^{-1}(t):=\inf\{s: \kappa(s)>t\}.
$$
Then it is easy to see that
$$
(\kappa^{-1})^{-1}=\kappa.
$$

Let $\sW^m:=C(\mR_+,\mR^m)$ be the space of continuous functions from $\mR_+$ to $\mR^m$,
endowed with compact uniform convergence topology.

We will need the following lemmas.

\bl\label{kurtz}
Let, for each $n$,  $y_n\in \sW^m$, $y\in \sW^m$
and $\t_n$ an increasing c\`adl\`ag function from $\mR_+$ to itself with $\t_n(0)=0$
and $\sup_n\t_n(t)<\infty$ for all $t$.
Suppose $y_n\longrightarrow y$ in $\sW^m$.
Define $x_n(t):=y_n(\t_n(t))$.
Then there exist a subsequence $\{n_k\}$ and an increasing c\`adl\`ag function $\t$
such that $\t_{n_k}$ converges to $\t$ at each continuity point of $\t$ and
$x_{n_k}(t)$ converges to $x(t)$ for all but countably many $t$ where $x(t):=y(\t(t))$.
If furthermore $y(u)=y(v)$ provided that $\t^{-1}(u)=\t^{-1}(v)$
$($or, equivalently, $y(\t(t))=y(\t(t-))$ for all $t$$)$, then $x\in \sW^m$
and $x_{n_k}\longrightarrow x$ in $\sW^m$.
\el

\pf Since $\sup_n\t_n(t)<\infty$ for all $t$, it follows from Lemma \ref{com} that $\{\t_n\}$ is compact in $\cT$. Consequently the
theorem follows directly from (the proof of) \cite[Lemma 2.3]{kur}.
(Notice that in that lemma $y_n$ are c\`adl\`ag functions and the condition ``strictly increasing" is
required, but it is clear from the proof that if we are restricted
to the continuous function space $\sW^m$, the strictness can be dropped).
\ef

\bl\label{inverse}
If $\t_n(t)\rightarrow \t(t)$ at each continuity point of $\t$,
and $\t$ and $\t_n$ are strictly increasing, then
$\t_n^{-1}\rightarrow \t^{-1}$ in $C(\mR_+,\mR_+)$.
\el

\pf
We first prove $\t_n^{-1}\rightarrow \t^{-1}$ pointwisely.
Let $\t^{-1}(t)=s$. By the strict increasingness of $\t$
we have $\t(u)>t$ for $u>s$ and $\t(u)<t$ for $u<s$.
Now for any $\e>0$, choose two points $s-\e<s_1<s<s_2<s+\e$ of continuity of $\t$.
Then for large $n$ we have
$$
\t_n(s_1)<t<\t_n(s_2)
$$
which implies
$$
s-\e\leq \t_n^{-1}(t)\leq s+\e.
$$
Now we prove the convergence is in fact uniform on each finite interval.
For simplicity we do this for $[0,1]$. For every $\e>0$, chose an $m$ such that
$$
\t^{-1}(k2^{-m})-\t^{-1}((k-1)2^{-m})<\e, ~~\forall k=1,2,\cdots, m.
$$
For large $n$ we have
$$
|\t^{-1}_n(k2^{-m})-\t^{-1}(k2^{-m})|<\e, ~~\forall  k=0,1,\cdots, m.
$$
Then, since $\t_n$ is increasing we have for all $t\in [0,1]$
$$
|\t^{-1}_n(t)-\t^{-1}(t)|\leq 6\e.
$$
This completes the proof.
\ef

Given a multivalued operator $A$ from $\mR^m$ to $2^{\mR^m}$, define:
\ce
D(A)&:=&\{x\in \mR^m:A(x)\neq \emptyset\}, \\
Im(A)&:=&\cup_{x\in D(A)}A(x),\\
Gr(A)&:=&\{(x,y)\in \mR^{2m}: x\in \mR^m, y\in A(x)\}.
\de
$A^{-1}$ is defined by: $x\in A^{-1}(y)\Leftrightarrow y\in A(x)$.

A maximal monotone operator $A$ is
a multivalued operaotr satisfying the following conditions:

(i) Monotonicity:
\ce
\<y_1-y_2,x_1-x_2\>\geq 0, \quad \forall (x_1,y_1),(x_2,y_2)\in Gr(A).
\de

(ii) Maximality:
\ce
(x_1,y_1)\in Gr(A) \Leftrightarrow\{\<y_1-y_2,x_1-x_2\>\geq 0,
\quad \forall (x_2,y_2)\in Gr(A)\}.
\de

Then we have ({see \cite{au}} or \cite{br})
\bp\label{Yosida}
(1)For each $x\in D(A)$, $A(x)$ is a closed and convex subset of $\mR^m$.
In particular, there is a unique $y\in A(x)$ such that $|y|$=$\inf\{|z|:z\in Ax\}$.
$A^\circ(x):=y$ is called the minimal section of $A$, and we have
\ce
x\in D(A)\Leftrightarrow |A^\circ(x)|<+\infty.
\de

(2) The resolvent operator $J_n:=(1+\frac{1}{n} A)^{-1}$ is
single-valued and Lipschitz continuous with Lipschitz constant 1.
Moreover, $\lim_{n\uparrow \infty}J_n x=x$ for any x $\in D(A)$.

(3) The Yosida approximation $A_n:=n(1-J_n)$ is monotone and
Lipschitz continuous with Lipschitz constant $n$. Moreover, as $n\uparrow \infty$
\ce
A_n(x)\rightarrow A^\circ(x)\quad and \quad |A_n(x)|\uparrow|A^\circ(x)|\quad if \quad x\in D(A)
\de
\ep

The following lemma which will be needed is proved in \cite{rx}.

\bl\label{barbu1}
If $x\notin D(A)$, $x_n\to x$, then
$$
\liminf_{n\to\infty}|A_n(x_n)|=\infty
$$
\el

We give the following definition for convenience.

\bd
Let $F,G$ be two continuous functions from $\mR_+$ to $\mR^m$
and suppose furthermore that $F(t)\in \overline{D(A)}$
for all $t$ and $G$ is of finite variation.
We say that $dG(t)\in A(F(t))dt$ if for every pair of continuous functions $(\alpha,\beta)$ satisfying
$$
(\alpha(t),\beta(t))\in Gr(A),
$$
we have
$$
\<F(t)-\alpha(t), dG(t)-\beta(t)dt\>\geq 0.
$$
\ed

Then we have the following result due to
\cite{ce2}.

\bl\label{fdg}
There exists $a\in \mR^m$, $c_1>0$, $c_2\geq 0$, such that if $dG(t)\in A(F(t))dt$ ,
then for all $0\leq s\leq t<\infty$, we have
$$
\int_s^t\<F(u)-a,dG(u)\>\geq c_1(|G|_t-|G|_s)-c_2\int_s^t|F(u)-a|du-c_1c_2(t-s).
$$
\el

Natural and important examples of maximal monotone operators
are subdifferentials of convex functions. More precisely,
Let $\vph$ be a proper convex function on $\mR^m$,
i.e., $\vph$ is a function from
$\mR^m$ to $(-\infty,+\infty]$ such that $\vph\not\equiv +\infty$ and
$$
\vph(\lambda x+(1-\lambda)y)\leq \lambda\vph(x)+(1-\lambda)\vph(y), ~~\forall \lambda\in (0,1).
$$
We also suppose that $\vph$ is lower-semicontinuus (l.s.c) and define its effective domain by
$$
D(\vph):=\{x\in \mR^m: \vph(x)<\infty\},
$$
and its subdifferential by
$$
\pvph(x):=\{y: \vph(x)\leq \vph(z)+(x-z,y)\}.
$$
We set
$$
D(\pvph):=\{x: \pvph(x)\neq\emptyset\}.
$$
Then it is well known that $\pvph$ is a multivalued maximal monotone operator and $D(\pvph)$ is a dense subset of $D(\vph)$ and $D(\vph)^o=D(\pvph)^o$ (see \cite{au, ba,br}).

For more examples of multivalued maximal monotone operators and applications of SVIs, we refer to  \cite{ce3}.

The following estimate, due to \cite{br}, will play a key role in this paper.

\bl\label{bis}
Suppose that $f\in L^2([0,T])$  and $u(t)$ be the unique solution to
the deterministic differential equation $\dot{u}(t)\in  -\pvph(u(t))+f$ with $u(0)\in D(\varphi)$.
Then we have
\begin{enumerate}[(i)]
\item $u$ is absolutely continuous in $[0,T]$;
\item
$$
\[\int_0^T|\dot{u}(t)|^2dt\]^{\frac{1}{2}}\leq \[\int_0^T|f(t)|^2dt\]^{\frac{1}{2}}+\sqrt{|\varphi(u(0))|}.
$$
\end{enumerate}
\el

\section{Main result}
To state our main result,
we have to fix some more notations first. Denote by
$$
\Om:=C_0([0,\infty),\mR^d)
$$
the space of continuous functions from $[0,\infty)$ to $\mR^d$ which are null at $0$,
and denote the generic point of $\Om$ by $\om$. For $\om\in\Om$, set
$$
w_t(\om):=\om(t).
$$
Endow $\Om$ with the compact uniform convergence topology. Denote
by $\cF$ the associated Borel $\sigma$-algebra and set
$$
\cF_t:=\sigma(w_s,s\leq t).
$$
Let $P$ stand for the canonical Wiener measure on $(\Om, \cF)$. Then $w_t$ is a standard
Brownian motion on  $(\Om, \cF, P)$.

Set
$$
\sS:=\{f\in \Om: f\mbox{ is smooth}\},
$$
$$
\sS_p:=\{f\in \Om: f\mbox{ is piecewise smooth}\}.
$$

From now on, $b$ and $\sigma$ will be $C^3_b$-maps from $\mR^m$ to $\mR^m$ and $\mR^m\times \mR^d$
respectively, $\vph$ is a proper l.s.c convex function on $\mR^m$ with $D(\pvph)^o\neq\emptyset$. We suppose further that
$D(\vph)$ is closed and $\vph$ is bounded on it (which in particular implies that $D(\vph)$
is closed). Consider the following Stratonovich SVI:
\be\label{eq10}
\left\{
\begin{array}{ll}
dX(t)\in b(X(t))dt+\sigma(X(t))\circ dw(t)-\pvph(X(t))dt, ~~t\in \mR_+\\
X(0)=x\in \overline{D(\pvph)}.
\end{array}
\right.
\ee
We have the following definition from \cite{ce, ce2}:
\bd\label{defsolution}
A pair of continuous and $\sF_t-$adapted
processes $(X,K)$ is called a solution of (1) if

(i) $X(0)=x$ and $X(t)\in \overline{D(\pvph)}\quad a.s.$;

(ii) $K=\{K(t), t\in \mR_+\}$ is of finite variation and $K(0)=0 \quad a.s.$;

(iii) $dX(t)=b(X(t))dt+\sigma(X(t))\circ dw(t)-dK(t)$, $t\in \mR_+$, \quad a.s.;

(iv) almost surely, $dK(t)\in \pvph(X(t))dt$.
\ed
By \cite{ce, ce2}, Eq.(\ref{eq10}) has a unique solution $(X(t),K(t))$.

For $h\in \sS_p$, consider the following deterministic variational inequality
\be\label{eqdeter}
\left\{
\begin{array}{ll}
d\xi(t)\in b(\xi(t))dt+\sigma(\xi(t))\dot{h}(t)dt-\pvph(\xi(t))dt, ~~t\in \mR_+,\\
\xi(0)=x\in \overline{D(\pvph)}.
\end{array}
\right.
\ee
By a classical result in \cite{br}, this inequality admits a unique solution and we shall denote it by
$\xi(h, x)$. Set then
$$
\eta(h,x,t):=\int_0^tb(\xi(h,x,s))ds+\int_0^t\sigma(\xi(h,x,s))\dot{h}(s)ds-\xi(h, x, t)
$$
$$
\sS^x:=\{(\xi(h,x), \eta(h,x)): h\in\sS\}.
$$
$$
\sS^x_p:=\{(\xi(h,x), \eta(h,x)): h\in\sS_p\}.
$$

We can now state our main result:

\bt\label{main}
Denote by $P_x$ the law of $(X, K)$, the unique solution of $($\ref{eq0}$)$, in $\sW^{2m}$. Then
$$
\mbox{supp}(P_x)=\overline{\sS_p^x}=\overline{\sS^x}.
$$
\et

This theorem will be a direct consequence of Theorems \ref{thsection3} and \ref{thsection4}.

\section{Limit theorem}

Now consider the following deterministic variational inequality
\be\label{apprsolution1}
\left\{
\begin{array}{ll}
\dot{X}_n(t)\in b(X_n(t))+\sigma(X_n(t))\dotw-\pvph(X_n(t)), \\
X_n(0)=x \in \overline{D(\pvph)},
\end{array}
\right.
\ee
where
$$
\dotw=2^n[w(t_n^+)-w(t_n)], \quad t_n^+=\frac{[2^nt]+1}{2^n}, \quad t_n=\frac{[2^nt]}{2^n}.
$$
Here $[a]$ stands for the integer part of $a$.

By \cite[Propsition 3.12]{br}, Eq.(\ref{apprsolution1}) has a unique solution  $X_n$. Set
$$
K_n(t)=\int_0^t b(X_n(s))ds+\int_0^t\sigma(X_n(s))\dot{w}_n(s) ds-X_n(t).
$$
Since $t\rightarrow X_n(t)$ is continuous and $x\in \overline{D(\pvph)}$, we have $X_n(t)\in D(\pvph)$ for $t\in \mR^m_+$.
Then we have by Lemma \ref{bis}
\be
\int_{u_n}^u|\dot{X}_n(t)|^2dt&\leq&
\int_{u_n}^u |b(X_n(v))+\sigma(X_n(v))\dot{w}_n(v)|^2dv+C\no\\
&\leq& C(1+2^{-n}|\dot{w}_n(u)|^2),
\ee
where $C$ is a constant independent of $n$.

Now we can state the main result of this section.
\bt
$(X_n,K_n)$ converges in $\sW^{2m}$ to $(X,K)$  in probability.
\et

The rest of this section is devoted to the proof of this theorem.
First we note that it is easy to deduce from Lemma \ref{com} the following tightness criterion (see also \cite{ce}):
\bl \label{tight1}
Let $\{\t_n\}_{n\in \mN^*}$ be a family of c\`adl\`ag increasing processes with $\t_n(0)=0$.
If for all $t>0$ there exists a constant $C(t)>0$ such that
$$
\sup_{n\in \mN^\ast}\mE[\t_n(t)]\leq C(t),
$$
then $(\t_n)_{n\in \mN^\ast}$ is tight on $\cT$.
\el

We set
$$\t_n(t):=|K_n|_t+t,$$
$$\tau_n:=\t_n^{-1},$$
$$M_n(t):=\int_0^t\sigma(X_n(u))\dotwu du,$$
$$Y_n(t):=X_n(\tau_n(t)),$$
$$H_n(t):=K_n(\tau_n(t)), $$
$$w_n(t):=\int_0^t\dotwu du,$$
$$a^{ij}(x):=\sum_{k=1}^d\sigma_k^i(x)\sigma_k^j(x),$$
$$\a_n(t):=\sigma(X_n(t_n))\dot{w}_n(t),$$
$$(\sigma'\sigma)_i^{l,l'}(x):=\sum_{j=1}^m\left(\p_j\sigma^{il}(x)\right)\sigma^{jl'}(x),$$
$$
(Lf)(x):=\sum_{i=1}^m \left[b_i(x)+\sum_{k=1}^d\sum_{j=1}^m(\frac{\p}{\p x^j}\sigma_k^i(x))\sigma_k^j(x)\right] \p_i f(x)
+\frac{1}{2}\sum_{i,j=1}^m a^{ij}(x)\p_i\p_jf(x).
$$
We now prove
\bt \label{tightness}
$(H_n, M_n, Y_n, w_n,\t_n)_{n\in \mN^\ast}$ is tight in $\sW^{3m} \times \Om\times \cT$.
\et
\begin{proof}Since
$$
t-s=\t_n(\tau_n(t))-\t_n(\tau_n(s))=|K_n|_{\tau_n(t)}-|K_n|_{\tau_n(s)}
+\tau_n(t)-\tau_n(s),
$$
we have for $s\leq t$
\be\label{tightness2}
|H_n(t)-H_n(s)|&=&|K_n(\tau_n(t))-K_n(\tau_n(s))|\no\\
&\leq &|K_n|_{\tau_n(t)}-|K_n|_{\tau_n(s)}\no\\
&=&t-s-(\tau_n(t)-\tau_n(s))\no\\
&\leq & t-s.
\ee
Thus the tightness of $\{H_n\}$ follows. That of $\{w_n\}$ is trivial since
\be
\mE[|w_n(t)-w_n(s)|^{2p}]\leq C|t-s|^p,~~\forall p\geq 1\label{tig}.
\ee

Next we look at $\t_n$. Let $a$ be as in Proposition \ref{fdg}.
We have for all $0\leq s\leq t<\infty$,
\ce
|X_n(s)-a|^{2}&=&|x-a|^{2}+\\
&&+2\int_0^s\<X_n(u)-a,b(X_n(u))\>du\\
&&+2\int_0^s\<X_n(u)-a,\sigma(X_n(u))\dotwu\>du\\
&&-2\int_0^s\<X_n(u)-a,dK_n(u)\>.
\de
By Proposition \ref{fdg}
\be\label{tightness4}
|X_n(s)-a|^{2}&\leq&C+Ct+C\int_0^s|X_n(u)-a|^{2}du\no\\
&&+C\int_0^s\<X_n(u)-a,\sigma(X_n(u))\dotwu\>du-C|K_n|_s,
\ee
where we have used the boundedness of $b$ and Young's inequality.
Let
\ce
I:=\mE\left[\sup_{0\leq s\leq t}\left|\int_0^s f(X_n(u))\dotwu du\right|\right],
\de
where $f(X_n(u))=(X_n(u)-a)^*\sigma(X_n(u))$. Then
\be
I&\leq &\mE\left[\sup_{0\leq s\leq t}\left|\int_0^sf(X_n(u_n))\dotwu du\right|\right]\no\\
&&+\mE\left[\sup_{0\leq s\leq t}\left|\int_0^s\int_{u_n}^u\frac{\p f(X_n(v))}{\p v}\dotwu dvdu\right|\right]\no\\
&:=&I_1+I_2.\label{idecompose}
\ee
Noticing that
$$
\int_0^sf(X_n(u_n))\dotwu du=\int_0^{s_n^+}\xi_udw(u),
$$
where
$$
\xi_u=2^n\int_{u_n}^{u_n^+\wedge s}f(X_n(v_n))dv=2^n(u_n^+\wedge s-u_n)f(X_n(u_n)).
$$
We have by BDG inequality and the boundedness of $\s$
\be\label{tightness5}
I_1&\leq & CE\left(\int_0^{t_n^+}|\xi_u|^2du\right)^\frac12\no\\
&\leq& C\mE\left( \int_0^t|X_n(u_n)-a|^{2}du \right)^\frac{1}{2}\no\\
&\leq& C\int_0^t\mE\[\sup_{0\leq u\leq s}|X_n(u)-a|^{2} \]ds+C(t)
\ee

By $\dot{X}_n(t)=b(X_n(t))+\sigma(X_n(t))\dot{w}_n(t)-\dot{K}_n(t)$
and the boundedness of $b$ and $\s$, we have
\be\label{tightness6}
I_2&\leq&\mE\Big[\sup_{0\leq s\leq t}\Big|\int_0^s\int_{u_n}^u
\left(|X_n(v)-a||\nabla\sigma(X_n(v))|+|\sigma(X_n(v))|\right)\no\\
&&\left(|b(X_n(v))|+|\sigma(X_n(v))||\dot{w}_n(v)|+|\dot{K}_n(v)|\right)|\dotwu|dvdu\Big|\Big]\no\\
&\leq&\mE\left[\sup_{0\leq s\leq t}\left|\int_0^s\int_{u_n}^u
\left(C|X_n(v)-a|+C\right)
\left(C+C|\dot{w}_n(v)|+|\dot{K}_n(v)|\right)|\dotwu|dvdu\right|\right]\no\\
&\leq&C\mE\left[\int_0^t\int_{u_n}^u |\dotwu|dvdu\right]
+C\mE\left[\int_0^t\int_{u_n}^u |\dotwu|^2dvdu\right]\no\\
&&+C\mE\left[\int_0^t\int_{u_n}^u |\dot{K}_n(v)||\dotwu|dvdu\right]
+C\mE\left[\int_0^t\int_{u_n}^u |X_n(v)-a||\dotwu|dvdu\right]\no\\
&&+C\mE\left[\int_0^t\int_{u_n}^u |X_n(v)-a||\dotwu|^2dvdu\right]\no\\
&&+C\mE\left[\int_0^t\int_{u_n}^u |X_n(v)-a||\dot{K}_n(v)||\dotwu|dvdu\right]\no\\
&\leq&C\mE\left[\int_0^t\int_{u_n}^u |\dotwu|dvdu\right]
+C\mE\left[\int_0^t\int_{u_n}^u |\dotwu|^2dvdu\right]\no\\
&&+C\mE\left[\int_0^t\int_{u_n}^u |\dot{X}_n(v)||\dotwu|dvdu\right]
+C\mE\left[\int_0^t\int_{u_n}^u |X_n(v)-a||\dotwu|dvdu\right]\no\\
&&+C\mE\left[\int_0^t\int_{u_n}^u |X_n(v)-a||\dotwu|^2dvdu\right]\no\\
&&+C\mE\left[\int_0^t\int_{u_n}^u |X_n(v)-a||\dot{X}_n(v)||\dotwu|dvdu\right]\no\\
&:=&I_{21}+I_{22}+I_{23}+I_{24}+I_{25}+I_{26}.
\ee
It is easily seen that
\be
I_{21}\leq C(t)2^{-\frac{n}{2}}
\ee
and
\be
I_{22}\leq C(t).
\ee
By Lemma \ref{bis} and the boundedness of $b$ and $\s$,
\be
I_{23}&\leq& C\mE\left[\int_0^t\left\{\int_{u_n}^u |\dot{X}_n(v)|^2dv\right\}^{\frac{1}{2}}
\left\{\int_{u_n}^u|\dotwu|dv\right\}^{\frac{1}{2}}du\right]\no\\
&\leq& C\mE\left[\int_0^t\left\{\int_{u_n}^u |b(X_n(v))+\sigma(X_n(v))\dot{w}_n(v)|^2dv\right\}^{\frac{1}{2}}
2^{-\frac{n}{2}}|\dot{w}_n(u)|du\right]+C(t)\no\\
&\leq& C\mE[|\dot{w}_n(t)|^22^{-n}]+C(t)\no\\
&\leq& C(t).
\ee
Moreover, we have
\be
I_{24}&=&C\int_0^t\int_{u_n}^u \mE[|X_n(v)-a||\dotwu|]dvdu\no\\
&\leq& C\int_0^t\int_{u_n}^u (\mE|X_n(v)-a|^2)^{\frac{1}{2}}
(\mE|\dotwu|^{2})^{\frac{1}{2}}dvdu\no\\
&\leq& C2^{-\frac{n}{2}}\int_0^t \mE\Big[\sup_{0\leq v\leq u}
|X_n(v)-a|^{2}\Big]du+C(t)2^{-\frac{n}{2}}
\ee
and
\be
I_{25}&=&C\int_0^t\int_{u_n}^u \mE[|X_n(v)-a||\dotwu|^2]dvdu\no\\
&\leq& C\int_0^t\int_{u_n}^u
(\mE[|X_n(v)-a|^{2}])^{\frac{1}{2}}
(\mE[|\dotwu|^{4}])^{\frac{1}{2}}dvdu\no\\
&\leq& C\int_0^t \mE\Big[\sup_{0\leq v\leq u}|X_n(v)-a|^{2}\Big]du+C(t)\label{tightness7}.
\ee
Furthermore,
\be
I_{26}&\leq& C\mE\left[\int_0^t\left\{\int_{u_n}^u |X_n(v)-a|^{2}dv\right\}^{\frac{1}{2}}
\left\{\int_{u_n}^u|\dot{X}_n(v)|^2|\dotwu|^2dv\right\}^{\frac{1}{2}}du\right]\no\\
&\leq& C\mE\left[\int_0^t\left\{\int_{u_n}^u |X_n(v)-a|^{2}dv\right\}^{\frac{1}{2}}|\dotwu|
\left\{\int_{u_n}^u|\dot{X}_n(v)|^2dv\right\}^{\frac{1}{2}}du\right]\no\\
&\leq& C2^{-\frac{n}{2}}\mE\left[\int_0^t\left\{\int_{u_n}^u |X_n(v)-a|^{2}dv\right\}^{\frac{1}{2}}|\dotwu|^2du\right]\no\\
&&+C2^{-\frac{n}{2}}\mE\left[\int_0^t\left\{\int_{u_n}^u |X_n(v)-a|^{2}dv\right\}^{\frac{1}{2}}|\dotwu|du\right]\no\\
&\leq&C\int_0^t \mE\Big[\sup_{0\leq v\leq u}|X_n(v)-a|^{2}\Big]du+C(t) \label{tightness8}
\ee

Combining (\ref{idecompose}), (\ref{tightness5}), (\ref{tightness6})-(\ref{tightness7}),
(\ref{tightness8}) gives
$$
I\leq C\int_0^t \mE\[\sup_{0\leq v\leq u}|X_n(v)-a|^{2}\]du+C(t).
$$
Hence
$$
\mE\[\sup_{0\leq s\leq t}|X_n(s)-a|^{2}\]\leq C(t)+
C\int_0^t \mE\[\sup_{0\leq v\leq u}|X_n(v)-a|^{2}\]du.
$$
Using Gronwall's inequality,
\be
\mE\[\sup_{0\leq s\leq t}|X_n(s)-a|^{2}\]\leq C(t)e^{Ct}. \label{bounded}
\ee
We obtain by (\ref{tightness4})
\be
\mE [|K_n|_t]\leq C(t)e^{Ct}. \label{abound}
\ee
That is
\be\label{tightness3}
\sup_{n\in \mN^\ast}\mE[\t_n(t)]\leq C(t)<\infty, \quad 0\leq t<\infty.
\ee
Therefore, in virtue of Lemma \ref{tight1}, $\t_n(t)$ is tight.

$\forall p\geq 1$, we also have
\ce
\mE[|M_n(t)-M_n(s)|^{2p}]&=&\mE\left[\left|\int_s^t\sigma(X_n(u))\dotwu du\right|^{2p}\right]\\
&\leq&C\mE\left[\left|\int_s^t\sigma(X_n(u_n))\dotwu du\right|^{2p}\right]\\
&&+C\mE\left[\left|\int_s^t(\sigma(X_n(u))-\sigma(X_n(u_n)))\dotwu du\right|^{2p}\right]\\
&\leq& C(t-s)^p+C\mE\left[\left|\int_s^t\int_{u_n}^u\frac{\p}{\p x_j}\sigma^{il}
(X_n(v))b^j(X_n(v))\dot{w}_n^l(u)dvdu\right|^{2p}\right]\\
&&+C\mE\left[\left|\int_s^t\int_{u_n}^u(\sigma'\sigma)^{l,l'}(X_n(v))\dot{w}_n^l(u) \dot{w}_n^{l'}(u)dvdu\right|^{2p}\right]\\
&&+C\mE\left[\left|\int_s^t\int_{u_n}^u\frac{\p}{\p x_j}\sigma^{il}(X_n(v))\dot{K}_n^j(v)\dot{w}_n^l(u)dvdu\right|^{2p}\right]\\
&\leq& C(t-s)^p+C(t-s)^{2p}\\
&&+C(t-s)^{2p-1}\mE\left[\int_s^tdu\left|\int_{u_n}^u|\dot{K}_n(v)||\dotwu|dv\right|^{2p}\right]\\
&\leq& C(t-s)^p+C(t-s)^{2p}\\
&&+C(t-s)^{2p-1}\mE\left[\int_s^tdu\left|\int_{u_n}^u|\dotwu|dv\right|^{2p}\right]\\
&&+C(t-s)^{2p-1}\mE\left[\int_s^tdu\left|\int_{u_n}^u|\dotwu|^2dv\right|^{2p}\right]\\
&&+C(t-s)^{2p-1}\mE\left[\int_s^tdu\left|\int_{u_n}^u|\dot{X}_n(v)||\dotwu|dv\right|^{2p}\right]\\
&\leq& C(t-s)^p+C(t-s)^{2p}\\
&&+C(t-s)^{2p-1}\mE\left[\int_s^t\left(\int_{u_n}^u|\dot{X}_n(v)|^2dv\right)^p\left(\int_{u_n}^u|\dot{w}_n(v)|^2dv\right)^pdu\right]\\
&\leq& C(t-s)^p+C(t-s)^{2p}\\
&&+Cn^{2p}2^{-2np}(t-s)^{2p-1}\mE\left[\int_s^t\left(\int_{u_n}^u|\dot{X}_n(v)|^2dv\right)^p2^{-np}|\dot{w}_n(v)|^{2p}du\right]\\
&\leq & C(t-s)^p+C(t-s)^{2p}.
\de
That is
\be
\mE[|M_n(t)-M_n(s)|^{2p}]\leq C(t-s)^p.
\label{mbound}
\ee
Furthermore,
\ce
\mE[|Y_n(t)-Y_n(s)|^{2p}]&=&\mE[|X_n(\tau_n(t))-X_n(\tau_n(s))|^{2p}]\\
&\leq&C\mE[|\int_{\tau_n(s)}^{\tau_n(t)}b(X_n(u))du|^{2p}]
+C\mE[|M_n(\tau_n(t))-M_n(\tau_n(s))|^{2p}]\\
&&+C\mE[|H_n(t)-H_n(s)|^{2p}]\\
&\leq&C\mE[|\tau_n(t)-\tau_n(s)|^{2p}]+C(t-s)^{2p}\\
&&+C\mE[|M_n(\tau_n(t))-M_n(\tau_n(s))|^{2p}].
\de
Using $|\tau_n(t)-\tau_n(s)|\leq |t-s|$, similarly to (\ref{mbound}), we
have
\be\label{est4}
\mE[|M_n(\tau_n(t))-M_n(\tau_n(s))|^{2p}]
&=&\mE\left[\left|\int_{\tau_n(s)}^{\tau_n(t)}\sigma(X_n(u))\dotwu du\right|^{2p}\right]\no\\
&\leq&C\mE\left[\left|\int_{\tau_n(s)}^{\tau_n(t)}\sigma(X_n(u_n))\dotwu du\right|^{2p}\right]\no\\
&&+C\mE\left[\left|\int_{\tau_n(s)}^{\tau_n(t)}(\sigma(X_n(u))-\sigma(X_n(u_n)))\dotwu du\right|^{2p}\right]\no\\
&\leq& C(t-s)^p.
\ee
Hence
\be\label{est501}
\mE[|Y_n(t)-Y_n(s)|^{2p}]\leq C(t-s)^p.
\ee
Combining (\ref{tightness2}), (\ref{tig}), (\ref{tightness3}), (\ref{mbound}),
(\ref{est501}) gives the desired tightness by Aldous's theorem (see \cite{ja}).
\end{proof}

Denote by $\{L_n,n\in N^*\}$ the distribution of
$(\tau_n, H_n, M_n, Y_n, w_n,\t_n)$ on $\sW^{3m+1}\times\Om\times\cT$.
Since $\sW^{3m+1}\times\Om\times \cT$ is a Polish space,
by Prokhorov's theorem and Proposition \ref{tightness},
there exists $L_{n_k}$ and Probability $L$ on $\sW^{3m+1}\times\Om\times \cT$
such that $L_{n_k}\rightarrow L( k\rightarrow \infty)$.
To simplify the notation, we suppose that $L_n\rightarrow L(n\rightarrow \infty)$.
By Skorohod's representation theorem, there exists a probability space
$(\hat{\Omega},\hat{\sF},\hat{P})$ on which are defined random variables
$(\hat{\tau}_n, \hat{H}_n, \hat{M}_n, \hat{Y}_n,\hat{w}_n,\hat{\t}_n)$,
$(\hat{\tau}, \hat{H}, \hat{M},\hat{ Y}, \hat{w},\hat{\t})$
such that
\be
(\hat{\tau}_n, \hat{H}_n, \hat{M}_n, \hat{Y}_n,\hat{w}_n,\hat{\t}_n)\sim (\tau_n,
H_n, M_n, Y_n,w_n,\t_n), \label{equiv}
\ee
$$
\hat{P}^{(\hat{\tau}, \hat{H}, \hat{M}, \hat{ Y}, \hat{w},\hat{\t})}=L,
$$
and as $n\rightarrow \infty$,
\be\label{wconverge}
(\hat{\tau}_n, \hat{H}_n, \hat{M}_n, \hat{Y}_n,\hat{w}_n,\hat{\t}_n)
\longrightarrow (\hat{\tau}, \hat{H}, \hat{M}, \hat{Y},\hat{w}, \hat{\t}) \quad a.s.
\ee
in $\sW^{3m+1}\times\Om\times \cT$.
Define
\be\label{define}
\hat{X}_n(t):=\hat{Y}_n(\hat{\t}_n(t)),\quad
\hat{X}(t):=\hat{Y}(\hat{\t}(t)),\quad
\hat{K}_n(t):=\hat{H}_n(\hat{\t}_n(t)), \quad
\hat{K}(t):=\hat{H}(\hat{\t}(t)).
\ee

We now pass from the convergence of $\hat{Y}_n$ to that of $\hat{X}_n$.
First, note that according to Lemma \ref{kurtz}, this will be
done if we prove the following

\bt \label{meyer}
There exists $\hat{\Omega}_0\in \hat{\sF}$,
$\hat{P}(\hat{\Omega}_0)=1$ such that for all $\hat{\omega}\in \hat{\Omega}_0$,
if there exist $0\leq s\leq t<\infty$ satisfying
$\hat{\tau}(s)(\hat{\omega})=\hat{\tau}(t)(\hat{\omega})$,
then
$\hat{Y}(s)(\hat{\omega})=\hat{Y}(t)(\hat{\omega})$.
\et

\begin{proof}
The proof is the same as that of \cite[Th.3.2]{rx}, except the Step (A) there.
But this is even easier here. In fact, since $\hat{X}_n(t,\om)\in \overline{D(\pvph)}$ ,
there exists an $\hat{N}_1$ with $\hat{P}(\hat{N}_1)=0$ such that
$\hat{X}(t,\om)=\lim_{n\mapsto \infty}\hat{X}_n(t,\om)\in \overline{D(\pvph)},~~~~\forall (t,\om)\in \mR_+\times \hat{N}^c_1, ~~ a.s.$ and thus Step (A) is done.
\end{proof}

By Lemma \ref{kurtz} and Theorem \ref{meyer}, there is a
subsequence $\{n_k\}$ such that
\be
(\hat{X}_{n_k},\hat{K}_{n_k})\longrightarrow
(\hat{X},\hat{K})~ \mbox{in}~\sW^{2m}   \quad
a.s.. \label{xkcon}
\ee
Consequently, the image measure of $(\hat{X}_{n_k},\hat{K}_{n_k})$
has a weak limit $\mu$ in $\sW^{2m}$.

Denote by $\mu_{n}$ the law of $(X_{n},K_{n})$.
Since $(\hat{X}_n,\hat{K}_n)$ and $({X}_n,{K}_n)$ are identically
distributed, $\mu$ is a weak limit of $\{\mu_n\}$ in $\sW^{2m}$.

Starting from the very beginning with an arbitrary subsequence and repeating the above
reasoning, we know that any subsequence has a weakly convergent sub-subsequence and
we thus arrive at the following result.

\bt\label{relativecompact}
$\{\mu_n\}$ is relatively compact in $\sW^{2m}$.
\et

Next we shall prove that the whole sequence $\{\mu_n\}$ converges weakly to a
unique limit  and we shall identify this limit.

Let
$$
\sV^m:=\{V: \mR_+\mapsto \mR^m, V(0)=0, V \mbox{ is continuous and
of finite variation on compacts}\}.
$$
By Theorem \ref{relativecompact}, $\{\mu_n\}$ has a weak limit.
Using the equivalence of weak solution and martingale problem,
we shall prove:

\bt\label{martingaleprob} Suppose that $\mu$ is a weak limit of
$\{\mu_n\}$. Let $(x(\cdot),v(\cdot))$ be coordinate processes on $\sW^m\times \sV^m$. Then, under $\mu$,
$$
f(x(t))-f(x(s))+\int_s^t\<\nabla f(x(u)),dv(u)\>- \int_s^t
Lf(x(u))du
$$
is a martingale for all $f\in C_b^2$.
\et
\pf
By a density argument, it suffices to prove that
$$
\mE^\mu\[F\cdot(f(x(t))-f(x(s)))\]= \mE^\mu\left[F\cdot \int_s^t
\left(Lf(x(u))du-\<\nabla f(x(u)),dv(u)\>\right)\right]
$$
for all $f\in C_0^{\infty}(\mR^m), 0\leq s<t$,
and bounded $\sB_s(\sW^m)\times \sB_s(\sV^m)$
measurable $F: \Omega\mapsto \mR$.
Clearly, it will suffice to do this when $s$ and $t$ have the form $k/2^{N}$
and $F$ is bounded continuous, and $\sB_s(\sW^m)\times \sB_s(\sV^m)$ measurable.
Observe that
\ce
&&\mE^{\mu_n}\left[F\cdot\left(f(x(t))-f(x(s))+
\int_s^t\<\nabla f(x(u)),dv(u)\>\right)\right]\\
&=&\mE^{{\mu}_n}\left[F\cdot\int_s^t\<\nabla
f(x(u)),b(x(u))\>du\right]+
\mE^P\left[F\cdot\int_s^t\<\nabla f(X_n(u)),{\a}_n(u)\>du\right]\\
&&+\mE^P\left[F\cdot\int_s^t\<\nabla f(X_n(u)),\dot{M}_n(u)-\a_n(u)\>du\right]\\
&&+\mE^P\left[\int_s^t\<\nabla f(X_n(u)),dK_n(u)\>\right]\\
&:=&J_{1,n}+J_{2,n}+J_{3,n}+J_{4,n}. \de Clearly, \be
J_{1,n}\rightarrow \mE^\mu\left[F\cdot\int_s^t\<\nabla
f(x(u)),b(x(u))\>du\right], \ee and \be
J_{4,n}\rightarrow&&\mE^P\left[\int_s^t\<\nabla f(X(u)),dK(u)\>\right]\no\\
=&&\mE^\mu\left[F\cdot\int_s^t\<\nabla f(x(u)),dv(u)\>\right]. \ee
We have to consider $J_{2,n}$ and $J_{3,n}$. First we prove
\be\label{j2} J_{2,n}\rightarrow \mE^{{\mu}}\left[F\cdot \int_s^t
L_u^0 f(x(u))du\right], \ee where
$$
L_u^0=\frac{1}{2}a^{ij}(x)(\p ^2/\p x_i\p x_j).
$$

Let $H(x)$ denote the Hessian matrix of $f$. Since
$$
\mE^P[\a_n(u)|\sB_{u_n}(\sW^m)\times \sB_{u_n}(\sV^m)]=0,
$$
we have
\ce
J_{2,n}&=&\mE^P\left[F\cdot \int_s^t\<\nabla f(X_n(u_n)),\a_n(u)\>du\right]\\
&&+\mE^P\left[F\cdot \int_s^t\<\nabla f(X_n(u))
-\nabla f(X_n(u_n)),\a_n(u)\>du\right]\\
&=&\mE^P\left[F\cdot \int_s^t\<\nabla f(X_n(u))
-\nabla f(X_n(u_n)),\a_n(u)\>du\right]\\
&=&\mE^P\left[F\cdot \int_s^tdu\int_{u_n}^udv\<\dot{M}_n(v),H(X_n(v))\a_n(u)\>\right]\\
&&+\mE^P\left[F\cdot \int_s^tdu\int_{u_n}^udv\<b(X_n(v)),H(X_n(v))\a_n(u)\>\right]\\
&&-\mE^P\left[F\cdot \int_s^tdu\int_{u_n}^u\<H(X_n(v))\a_n(u),dK_n(v)\>\right]\\
&:=&K_{1,n}+K_{2,n}+K_{3,n}.
\de
An elementary calculus gives $|K_{2,n}|\rightarrow 0$.
Since
\ce
K_{3,n}\leq\mE\left[\left|\int_s^tdu\int_{u_n}^u\<H(X_n(v))\a_n(u),dK_n(v)\>\right|\right],\\
\de
we have to show that $\left|\int_s^tdu\int_{u_n}^u\<H(X_n(v))\a_n(u),dK_n(v)\>\right|$ is uniformly integrable
and converges to zero in probability.

Since
\ce
&&\int_{u_n}^u\<H(X_n(v))\a_n(u),dK_n(v)\>\\
&&=\<H(X_n(u))\a_n(u),K_n(u)-K_n(u_n)\>\\
&&-\int_{u_n}^u \left<\sum_{i=1}^m\a_n^i(u)\(\sum_{j=1}^m\frac{\p H_{ki}(X_n(v))}{x_j}\dot{X}_n(v)\),K_n(v)-K_n(u_n)\right>dv\\
&\leq& C|\dot{w}_n(u)||K_n(u)-K_n(u_n)|+\int_{u_n}^u|\dot{w}_n(u)||K_n(v)-K_n(u_n)||\dot{X}_n(v)|dv
\de
and
\ce
|K_n(u)-K_n(u_n)|&\leq&\int_{u_n}^u |\dot{K}_n(v)|dv\\
&\leq&\int_{u_n}^u |\dot{X}_n(v)|dv+C2^{-n}+C2^{-n}|\dot{w}_n(u)|\\
&\leq&2^{-\frac{n}{2}}\left\{\int_{u_n}^u |\dot{X}_n(v)|^2dv\right\}^{\frac{1}{2}}+C2^{-n}+C2^{-n}|\dot{w}_n(u)|\\
&\leq&C2^{-n}+C2^{-n}|\dot{w}_n(u)|,
\de
we have
\ce
&&\left|\int_{u_n}^u\<H(X_n(v))\a_n(u),dK_n(v)\>\right|\\
&\leq& C2^{-n}|\dot{w}_n(u)|+C2^{-n}|\dot{w}_n(u)|^2\\
&&+C2^{-n}\int_{u_n}^u|\dot{w}_n(u)||\dot{X}_n(v)|dv
+C2^{-n}\int_{u_n}^u|\dot{w}_n(u)|^2|\dot{X}_n(v)|dv\\
&\leq& C2^{-n}|\dot{w}_n(u)|+C2^{-n}|\dot{w}_n(u)|^2
+C2^{-2n}|\dot{w}_n(u)|\\
&&+C2^{-2n}|\dot{w}_n(u)|^2+C2^{-2n}|\dot{w}_n(u)|^3.
\de
For any $p>1$,
\ce
&&\mE\left[\left|\int_s^tdu\int_{u_n}^u\<H(X_n(v))\a_n(u),dK_n(v)\>\right|^p\right]\\
&\leq&C(t-s)^{p-1}\mE\left[\int_s^t\left|\int_{u_n}^u\<H(X_n(v))\a_n(u),dK_n(v)\>\right|^pdu\right]\\
&\leq& C(t-s)^p \{\mE[2^{-np}|\dot{w}_n(u)|^p]+\mE[2^{-np}|\dot{w}_n(u)|^{2p}]
+\mE[2^{-2np}|\dot{w}_n(u)|^p]\\
&&+\mE[2^{-2np}|\dot{w}_n(u)|^{2p}]+\mE[2^{-2np}|\dot{w}_n(u)|^{3p}]\}\\
&\leq&C.
\de
Therefore, $\left|\int_s^tdu\int_{u_n}^u\<H(X_n(v))\a_n(u),dK_n(v)\>\right|$ is uniformly integrable.

Set
$$
J_n(u):=\int_{u_n}^u \<H(X_n(v), dK_n(v)\>,
$$
$$
\beta_n(u):=\int_{u_n}^u \a_n(v)dv,
$$
$$
\gamma_n(u):=\int_{0}^u \a_n(v)dv.
$$
As $J_n(u_n)=0$, we have
\ce
&&\int_s^tdu\int_{u_n}^u\<H(X_n(v))\a_n(u),dK_n(v)\>\\
&=&\int_s^t\a_n(u)J_n(u)du\\
&=&\int_s^t J_n(u)d\gamma_n(u)\\
&=&\int_s^t J_n(u)d\beta_n(u)\\
&=&J_n(t)\beta_n(t)-J_n(s)\beta_n(s)-\int_s^t \beta_n(u)dJ_n(u).
\de
Since as $n\to\infty$
$$
|\beta_n(u)|\leq C|w(u_n^+)-w(u_n)|\rightarrow 0,
$$
$$
\sup_nE[V_t(J_n)]\leq C\sup_n E[V_t(K_n)]<\infty,
$$
we have
$$
\int_s^tdu\int_{u_n}^u\<H(X_n(v))\a_n(u),dK_n(v)\>\rightarrow 0
$$
in probability.
Hence $|K_{3,n}|\rightarrow 0$.

For $K_{1,n}$, we have
\ce
K_{1,n}&=&\mE^P\left[F\cdot \int_s^tdu\int_{u_n}^udv\<\a_n(v),H(X_n(v))\a_n(u)\>\right]\\
&&+\mE^P\left[F\cdot \int_s^tdu\int_{u_n}^udv\int_{v_n}^vdr
\<\frac{\p}{\p x_j}\sigma^{il}(X_n(r))b^j(X_n(r))\dot{w}_n^l,H(X_n(v))\a_n(u)\>\right]\\
&&+\mE^P\left[F\cdot \int_s^tdu\int_{u_n}^udv\int_{v_n}^vdr
\<(\sigma'\sigma)^{l,l'}(X_n(r))\dot{w}_n^l\dot{w}_n^{l'},H(X_n(v))\a_n(u)\>\right]\\
&&-\mE^P\left[F\cdot \int_s^tdu\int_{u_n}^udv\int_{v_n}^vdr
\<\frac{\p}{\p x_j}\sigma^{il}(X_n(r))\dot{K}_n^j(v)\dot{w}_n^l,H(X_n(v)) \a_n(u)\>\right]\\
&:=&K_{4,n}+K_{5,n}+K_{6,n}+K_{7,n}.
\de
Again, we have $|K_{5,n}|\rightarrow0$,
$|K_{6,n}|\rightarrow0$,$|K_{7,n}|\rightarrow0$ by simple calculus.
We only need to consider $K_{4,n}$.
\ce
K_{4,n}&=&\mE^P\left[F\cdot \int_s^tdu\int_{u_n}^udv\<\a_n(v),H(X_n(u_n)\a_n(u)\>\right]\\
&&+\mE^P\left[F\cdot \int_s^tdu\int_{u_n}^udv\<\a_n(v),(H({X}_n(v))-H(X_n(u_n)) \a_n(u)\>\right]\\
&:=&K_{8,n}+K_{9,n}.
\de
Since $|K_{9,n}|\rightarrow0$, it remains to examine $K_{8,n}$.
\ce
K_{8,n}&=&2^n\mE^P\left[F\cdot \int_s^tdu\int_{u_n}^udvtr[\sigma^\ast(X_n(v_n))
H(X_n(u_n))\sigma(X_n(u_n))]\right]\\
&=&2^n\mE^{{\mu}_n}\left[F\cdot \int_s^tdu\int_{u_n}^udvtr[\sigma^\ast(x(v_n))H(x(u_n)) \sigma(x(u_n))]\right].
\de
Since for bounded measurable functions $\vph$ and $\psi$  on $[s,t]$,
$$
2^n\int_s^t\vph(u)du\int_{u_n}^u\psi(v)dv \rightarrow \frac{1}{2}\int_s^t\vph(u)\psi(u)du
$$
(see \cite[Lemma 4.2]{st}), using ${\mu}_n\rightarrow{\mu}$,
$$
K_{8,n}\rightarrow \frac{1}{2}\mE^{{\mu}}\left[F\cdot \int_s^t tr[\sigma^\ast(x(u))H(x(u))
\sigma(x(u))]du\right],
$$
(\ref{j2}) is proved.

Finally we treat $J_{3,n}$. We write \ce
J_{3,n}&=&\mE^P\left[F\cdot\int_s^t\int_{u_n}^udv\<\nabla f(X_n(u)),
(\sigma'\sigma)^{l,l'}(X_n(v)) \dot{w}_n^l\dot{w}_n^{l'}\>\right]\\
&&+\mE^P\left[F\cdot\int_s^t\int_{u_n}^udv\<\nabla f(X_n(u)),
\frac{\p}{\p x_j}\sigma^{il}(X_n(v))b^j(X_n(v))\dot{w}_n^l\>\right]\\
&&-\mE^P\left[F\cdot\int_s^t\int_{u_n}^udv\<\nabla f(X_n(u)),
\frac{\p}{\p x_j}\sigma^{il}(X_n(v))\dot{K}_n^j(v)\dot{w}_n^l\>\right]\\
&:=&L_{1,n}+L_{2,n}+L_{3,n}. \de Clearly, $|L_{2,n}|\rightarrow 0$
and $|L_{3,n}|\rightarrow 0$. Observe that \ce
L_{1,n}&=&\mE^P\left[F\cdot\int_s^t\<\nabla f(X_n(u_n)),
(\sigma'\sigma)^{l,l'}(X_n(u_n))(u-u_n) \dot{w}_n^l\dot{w}_n^{l'}\>\right]\\
&&+\mE^P\left[F\cdot\int_s^tdu\int_{u_n}^udv\<\nabla f(X_n(u_n)),
[(\sigma'\sigma)^{l,l'}(X_n(v))-(\sigma'\sigma)^{l,l'}(X_n(u_n))]\dot{w}_n^l\dot{w}_n^{l'}\>\right]\\
&&+\mE^P\left[F\cdot\int_s^tdu\int_{u_n}^udv\<\nabla f(X_n(u))-
\nabla f(X_n(u_n)),(\sigma'\sigma)^{l,l'}(X_n(v)) \dot{w}_n^l\dot{w}_n^{l'}\>\right]\\
&:=&L_{4,n}+L_{5,n}+L_{6,n}. \de By $|L_{5,n}|\rightarrow 0$,
$|L_{6,n}|\rightarrow 0$ and $L_{4,n}\rightarrow
\frac{1}{2}\mE^{{\mu}} \left[F\cdot \int_s^t \< \nabla
f(x(u)),\sigma' \sigma(x(u)) \>du\right]$, we get \be\label{j3}
J_{3,n} \rightarrow \frac{1}{2}\mE^{{\mu}} \left[F\cdot \int_s^t \<
\nabla f(x(u)),\sigma' \sigma(x(u)) \>du\right]. \ee The proof is
completed.
\end{proof}

Using the same argument as in the proof of \cite[Prop. 5.13]{ce}, we can prove:
\bp\label{martingaleprob2}
If $(\a,\b)$ are continuous functions satisfying
$$(\a(t),\b(t))\in Gr(A),\quad \forall t\in \mR_+,$$
then the measure
$$\<x(t)-\a(t), dv(t)-\b(t)dt\>$$
is positive on $\mR^+$, $\mu$-a.s..
\ep

Now, instead of (\ref{eq10}) we consider the following system
\be
\begin{cases}
& dX^i(t)\in\sigma^{ij}(X(t))\circ dw^j(t)+b^i(X(t))dt-(\pvph)^i(X(t))dt, ~i=1,2,\cdots m\\
& dX^{m+j}(t)=dw^j(t)\\
& X^i(0)=x^i\\
& X^{m+j}(0)=0
\end{cases}
\ee

Denote by $\nu_n$ the law of $(w_n, X_n, K_n)$ in $\Om\times \sW^m\times \sV^m$.

Applying Theorem \ref{martingaleprob} and Propsition \ref{martingaleprob2} to
the above system and using the uniqueness in distribution of the solution of (\ref{eq10}) we obtain

\bt \label{weak} On $(\Om\times \sW^m\times \sV^m, \nu)$,
$t\mapsto w(t)$ is a Brownian motion and $(X,K)$ is
solution of the following multivalued Stratonovich SDE:
\be \left\{
\begin{array}{ll}
d X(t) \in b(X(t))dt+\sigma(X(t))\circ d w(t)-\pvph(X(t))dt, \\
X(0)=x \in \overline{D(\pvph)},
\end{array}
\right.
\ee
Moreover, $\nu$ is the unique weak limit of  $\{\nu_n\}$.
\et

Finally, using an argument analogous to \cite[Th.6.2]{rx} we can prove
\bt\label{thsection3}
$(X_n,K_n)$ converges in $\sW^{2m}$ to $(X,K)$  in probability.
\et

\section{approximate continuity}
In this section we further suppose that $\vph$ is bounded and the notations $a_{ij}(x)$,
$(\sigma'\sigma)_i^{l,l'}(x)$ and $(Lf)(x)$ which will be needed are defined as in Section 3.

Set
$$
H_0:=\mW^d_0\cap C_b^2.
$$
Let $h\in H_0$ and denote by $\xi(t)$ the unique solution, whose existence and uniqueness is assured by
\cite[Proposition 3.12]{br}, of the following DVI:
\be\label{mode2}
\left\{
\begin{array}{ll}
\dot{\xi}(t) \in b(\xi(t))+\sigma(\xi(t))\dot{h}(t) -\pvph(\xi(t)), \\
\xi(0)=x \in \overline{D(\pvph)}.
\end{array}
\right. \label{mode2}
\ee
Then
$$
\dot{\eta}(t):=-\dot{\xi}(t)+ b(\xi(t))+\sigma(\xi(t))\dot{h}(t) \in\pvph(\xi(t)).
$$

The following two Lemmas and corollary are taken from \cite{iw}.

\bl\label{w}
$P(\|w\|_T<\e)\sim C\exp(-\frac{C}{\e^2})$ as $\e\downarrow 0$.
\el

\bl
Set $\kappa^{ij}(t):=\frac{1}{2}\int_0^t[w^i(s)dw^j(s)-w^j(s)dw^i(s)]$, $i,j=1,...,d$.
Then for all $i,j=1,...,d$,
$$\lim_{M\uparrow \infty}\sup_{0<\d\leq 1}P(\|\kappa^{ij}\|_T>M\d|\|w\|_T<\d)=0.$$
\el

\bc \label{co}
Let $\zeta^{ij}(t):=\int_0^tw^i(s)\circ dw^j(s)$, $i,j=1,...,d$. Then for all $i,j=1,...,d$, we have
$$\lim_{M\uparrow \infty}\sup_{0<\d\leq 1}P(\|\zeta^{ij}\|_T>M\d|\|w\|_T<\d)=0.$$
In particular, for every $\e>0$ and $\alpha\in (0,1)$,
\be
P(\|\zeta^{ij}\|_T>\e\delta^\alpha|\|w\|_T<\d)\rightarrow 0 \quad as \quad \d\downarrow 0.  \label{est2}
\ee
\ec

We will need two more lemmas.
\bl\label{expint}
Let $|K|_T$ be the total variation of $K$ on $[0,T]$.
Then there exists a strictly positive constant $\a>0$, such that
\be
\mE\[e^{\a|K|_T^2}\]<\infty.  \label{exp}
\ee
\el
\begin{proof}
By Lemma \ref{fdg} we have
\ce
|X(t)-a|^2&\leq & C+Ct+2\int_0^t\<X(s)-a, \sigma(X(s))\circ d w(s)\>
+C\int_0^t |X(s)-a|ds\\
&&-2\int_0^t\<X(s)-a, dK(s)\>\\
&\leq& C+Ct+C\int_0^t |X(s)-a|ds-C|K|_t+2\int_0^t\<X(s)-a, \sigma(X(s))\circ d w(s)\>.
\de
Hence
\be
|K|_T\leq CT+2\sup_{0\leq t\leq T}|\int_0^t\<X(s)-a, \sigma(X(s))\circ d w(s)\>|\label{est0}
\ee
Let $$
N_t:=\int_0^t\<X(s)-a, \sigma(X(s))\circ d w(s)\>,
$$
Since $X$ and $\sigma$ are bounded, we have
$$
\<N,N\>_t\leq Ct,
$$
Hence there exists $\a_1>0$, such that
\be
\mE\left[e^{\a_2 \|N\|_T^2}\right]< \infty. \label{exp3}
\ee
Combining (\ref{est0}) and (\ref{exp3}) gives the desired result.
\end{proof}

\bl \label{lem0}
Let $M_t:=\int_0^t\<X(s)-a, \sigma(X(s))\circ dw(s)\>$, $t\in [0,T]$. Then for all $\e>0$, we have
\be
P\left(\|M\|_T\geq \e\d^{-\frac{1}{2}}|\|w\|_T<\d\right)\rightarrow 0, ~~~\d\downarrow 0. \label{est31}
\ee
\be
P\left(|K|_T\geq \e\d^{-\frac{1}{2}}|\|w\|_T<\d\right)\rightarrow 0, ~~~\d\downarrow 0. \label{est9}
\ee
\el
\begin{proof}
Obviously
\be
P\left(|K|_T\geq \e \d^{-\frac{3}{2}}|\|w\|_T<\d\right)\sim\frac{Ce^{-C\e \d^{-3}}}{Ce^{-C\d ^{-2}}}\rightarrow 0, ~~\d\downarrow0. \label{est5}
\ee

Applying Ito's formula to $M_t$, we have
\be\label{m1}
M_t&=&\sum_{i=1}^m\int_0^t(X^i(s)-a)\sigma_k^i(X(s))\circ dw^k(s)\no\\
&=&\sum_{i=1}^m (X^i(t)-a)\sigma_k^i(X(t))w^k(t)-\sum_{i=1}^m\int_0^tw^k(s)(X^i(s)-a)\circ d\sigma_k^i(X(s))\no\\
&&-\sum_{i=1}^m\int_0^tw^k(s)\sigma_k^i(X(s))\circ d(X^i(s)-a)\no\\
&=&\sum_{i=1}^m (X^i(t)-a)\sigma_k^i(X(t))w^k(t)\no\\
&&-\sum_{i=1}^m\int_0^tw^k(s)(X^i(s)-a)\sigma_{k,l}^i(X(s))
[\sigma_{\a}^l(X(s))\circ dw^{\a}(s)+b^l(X(s))ds-dK^l(s)]\no\\
&&-\sum_{i=1}^m\int_0^tw^k(s)\sigma_k^i(X(s))[\sigma_{\a}^i(X(s))\circ dw^{\a}(s)+b^i(X(s))ds-dK^i(s)]\no\\
&=:&\sum_{j=1}^7J_j.
\ee
We have to prove that $P(\|J_j\|_T\geq \e \d^{-\frac{1}{2}}|\|w\|_T< \d)\rightarrow 0$
as $\d\downarrow 0$, for $j=1,2,..7$.

Obviously, there is no problem for $J_1$, $J_3$ and $J_6$.
For $J_2$, we have
\ce
J_2&=&\sum_{i=1}^m\int_0^tw^k(s)(X^i(s)-a)\sigma_{k,l}^i(X(s))\sigma_{\a}^l(X(s))\circ dw^{\a}(s)\\
&=&\sum_{i=1}^m\int_0^tw^k(s)(X^i(s)-a)\sigma_{k,l}^i(X(s))\sigma_{\a}^l(X(s))dw^{\a}(s)\\
&&+\frac{1}{2}\sum_{i=1}^m\int_0^t [\sigma_{k,l}^i\sigma_{\a}^l\sigma_{\b}^i](X(s))w^k(s)\d ^{\a \b} ds\\
&&+\frac{1}{2}\sum_{i=1}^m\int_0^t(X^i(s)-a) \frac{\p}{\p x^n}[\sigma_{k,l}^i\sigma_{\a}^l](X(s))\sigma_{\b}^n(X(s))w^k(s)\d ^{\a\b} ds\\
&:=&J_{21}+J_{22}+J_{23}.
\de
It is obvious that $J_{21}$ is a martingle with
\ce
<J_{21},J_{21}>_t=\sum_{i=1}^m\sum_{i'=1}^m\int_0^t(X^i(s)-a)(X^{i'}(s)-a) [\sigma_{k,l}^i\sigma_{k',l'}^{i'} a^{ll'}](X(s))w^k(s)w^{k'}(s)ds.
\de
Therefore, if $\|w\|_T<\d$ then $<J_{21},J_{21}>\leq C\d^2.$ Hence
\ce
P(\|J_{21}\|_T>\e\d^{-\frac{1}{2}},\|w\|_T<\d)&\leq& P(\max_{0\leq t\leq C\d^2}|B(t)|>\e\d^{-\frac{1}{2}})\\
&=& P(\max_{0\leq t\leq 1}|B(t)|>C\e\d^{-\frac{3}{2}})\\
&=&2\int_{C\e\d^{-\frac{3}{2}}}^{\infty}\frac{1}{\sqrt{2\pi}}\exp\{-\frac{x^2}{2}\}dx\\
&\leq& C\exp\{-C\frac{\e^2}{\d^3}\}.
\de
By Lemma \ref{w}, $P(\|w\|_T\leq \d)\sim C\exp\{-C\frac{1}{\d^2}\}$. Hence
\ce
P(\|J_{21}\|_T>\e\d^{-\frac{1}{2}}|\|w\|_T<\d)\rightarrow 0.
\de
Since $P(\|J_{22}\|_T>\e\d^{-\frac{1}{2}}|\|w\|_T<\d)\rightarrow 0$ and
$P(\|J_{23}\|_T>\e\d^{-\frac{1}{2}}|\|w\|_T<\d)\rightarrow 0$ as $\d\downarrow 0$, we have
\ce
P(\|J_2\|_T>\e\d^{-\frac{1}{2}}|\|w\|_T<\d)\rightarrow 0.
\de
Similarly, we have
\ce
P(\|J_{5}\|_T>\e\d^{-\frac{1}{2}}|\|w\|_T<\d)\rightarrow 0.
\de
Since $|J_4|\leq C\|w\|_T|K|_T$ and $|J_7|\leq C\|w\|_T|K|_T$, we obtain  by (\ref{est5})
\ce
P(\|J_4\|_T>\e\d^{-\frac{1}{2}}|\|w\|_T<\d)\rightarrow 0
\de
and
\ce
P(\|J_7\|_T>\e\d^{-\frac{1}{2}}|\|w\|_T<\d)\rightarrow 0.
\de
Summing up gives
\be
P\left(\|M\|_T\geq \e \d^{-\frac{1}{2}}|\|w\|_T<\d\right)\rightarrow 0, ~~\d\downarrow0.  \label{est6}
\ee
Combining (\ref{est6}) with
\be\label{kbym}
|K|_T&\leq& CT+2\sup_{0\leq t\leq T}|\int_0^t\<X(s)-a, \sigma(X(s))\circ d w(s)\>|
\ee
yields
\be\label{k3on2}
P\left(|K|_T\geq \e \d^{-\frac{1}{2}}|\|w\|_T<\d\right)\rightarrow 0, ~~\d\downarrow0. \label{est7}
\ee
This completes the proof.
\end{proof}

Now we have
\bl\label{lemforevery}
For every $\e>0$,
$$
\lim_{\d\downarrow 0}P(\|\zeta^{km}\|_T|K|_T>\e|\|w\|_T<\d)=0.
$$
\el
\begin{proof}
\ce
&&P(\|\zeta^{km}\|_T|K|_T>\e|\|w\|_T<\d)\\
&\leq&P(\|\zeta^{km}\|_T>\d^{\frac{1}{2}}|\|w\|_T<\d)+P(|K|_T>\e \d^{-\frac{1}{2}}|\|w\|_T<\d)\\
& \rightarrow &0.
\de
\end{proof}

It follows immediately:
\bl \label{lem2}
\be
P\left(\left\|\int_0^t\zeta^{km}(s)dK(s)\right\|_T>\e|\|w\|_T<\d\right)\rightarrow 0. \label{est10}
\ee
\el

We also need the following
\bl \label{lem}
Let $f(x):\mR^m\rightarrow \mR$ be bounded and uniformly continuous. Then for all $\e>0$ and $k,m=1,2,...,d$
\be\label{es}
P\left(\left\|\int_0^tf(X(s))d\zeta^{km}(s)\right\|_T>\e|\|w\|_T<\d\right)\rightarrow 0 \quad as \quad \d\downarrow 0.
\ee
\el
\begin{proof}
First we assume that $f\in C_b^2(\mR^m)$. Set $f_l(x):=\frac{\p}{\p x_l}f(x)$. By It\^o's formula we have
\ce
\int_0^tf(X(s))d\zeta^{km}(s)&=&f(X(t))\zeta^{km}(t)
-\int_0^t\zeta^{km}(s)f_l(X(s))\sigma_i^l(X(s))dw^i(s)\\
&&-\int_0^t(Lf)(X(s))\zeta^{km}(s)ds-\int_0^tf_l(X(s))\sigma_m^l(X(s))w^k(s)ds\\
&&+\int_0^tf_l(X(s))\zeta^{km}(s)dK^l(s)\\
&:=&P_1+P_2+P_3+P_4+P_5.
\de
We now prove that
\be
\lim_{\d\downarrow 0}P(\|P_i\|_T>\e|\|w\|_T<\d)=0
\ee
for $i=1,\cdots, 5$. This is deduced plainly from (\ref{est2}) for $i=1,i=3$ and is trivial for $i=4$.
Since
$$
\left|\int_0^tf_l(X(s))\zeta^{km}(s)dK^l(s)\right|\leq C\|\zeta^{km}\|_T\cdot|K|_T,
$$
we obtain by Lemma \ref{lemforevery}
\be
\lim_{\d\downarrow0}P(\|P_5\|_T>\e| \|w\|_T<\d)=0.
\ee
It remains to look at  $P_2$.
Set $\a_i(x):=-f_l(x)\sigma_i^l(x)$, $\a_{i,l}:=\frac{\p}{\p x^l}\a_i(x)$. Then by It\^o's formula,
\ce
P_2&=&\int_0^t\a_i(X(s))\zeta^{km}(s)dw^i(s)\\
&=&\a_i(X(t))\zeta^{km}(t)w^i(t)-\int_0^t\a_{i,l}(X(s))\sigma_j^l(X(s))\zeta^{km}(s)w^i(s)dw^j(s)\\
&&-\int_0^t(L\a_i)(X(s))\zeta^{km}(s)w^i(s)ds
-\int_0^t\a_i(X(s))w^i(s)d\zeta^{km}(s)\\
&&-\int_0^t\a_m(X(s))w^k(s)ds-\int_0^t\zeta^{km}(s)\a_{i,l}(X(s))\sigma_j^l(X(s))\d^{ij}ds\\
&&-\int_0^t\a_{i,l}(X(s))\sigma_m^l(X(s))w^i(s)w^k(s)ds
+\int_0^t\a_{i,l}(X(s))\zeta^{km}(s)w^i(s)dK^l(s)\\
&:=&\sum_{i=1}^8L_i.
\de
Again it is sufficient to show that $P(\|L_i\|_T>\e| \|w\|_T<\d)\rightarrow 0$.
The proof can be done in a similar way to the estimation of $I_2$ in \cite[p. 522-524]{iw} except an extra
term $L_8$. But $L_8$ can be estimated easily by  using Lemma \ref{lemforevery}.

The passage from $C_b^2$ functions to bounded and uniformly continuous functions is
completely the same as on \cite[p.525]{iw} and so we omit it.
\end{proof}

\bl\label{fbecb2}
Let $f$ be a $C_b^2$-function. Then for $k=1,\cdots, d$,
$$
\lim_{\d\downarrow 0}P(\|\int_0^tf(X(s))\circ dw^k(s)\|_T>\e|\|w\|_T<\d)=0.
$$
\el
\begin{proof}
Set $f_l:=\partial_l f$. By Ito's formula we have
\ce
&&\int_0^tf(X(s))\circ dw^k(s)\\
&=& f(X(t))\circ dw^k(t)-\int_0^t[f_l\sigma_m^l](X(s))w^k(s)\circ dw^m(s)\\
&&-\int_0^t[f_lb^l](X(s))w^k(s)ds-\int_0^tf_l(X(s))w^k(s)dK^l(s)\\
&=&I_1(t)+I_2(t)+I_3(t)+I_4(t).
\de
It is sufficient to prove
$$
\lim_{\d\downarrow 0}P(\|I_i(t)\|_T>\e|\|w\|_T<\d)=0, ~\forall i=1,2,3,4.
$$
$I_1$ and $I_3$  are trivial and $I_4$ will make no trouble by Lemma \ref{lem0}. It remains to look at $I_2$. We have
\ce
I_2(t)&=&-\int_0^t[f_l\sigma_m^l](X(s))\circ d\zeta^{km}(s)\\
&=&-\int_0^t[f_l\sigma_m^l](X(s))d\zeta^{km}(s)-\frac12\int_0^t\frac{\partial}{\partial x_j}[f_l\sigma_m^l]\sigma_q^j(X(s))w^k(s)\d^{qm} ds\\
&=&J_1(t)+J_2(t)
\de
We claim that
$$
\lim_{\d\downarrow 0}P(\|J_i(t)\|_T>\e|\|w\|_T<\d)=0, ~\forall i=1,2.
$$
For $i=2$ this is obvious and for $i=1$ this follows from Lemma \ref{lem}. The proof is now complete.
\end{proof}
Now we are in the position to state our main result of this section.
\bt\label{thsection4}
$\forall h\in \sS$ and $\e>0$
\be
P(\|X(t)-\xi(t)\|_T+|K(t)-\eta(t)|_T<\e|\|w-h\|_T<\d)\rightarrow 1 \quad as \quad \d\downarrow 0. \label{conti}
\ee
\et

\begin{proof}
By the standard argument (see \cite{st} or \cite[p.527-528]{iw}) it suffices to prove (\ref{conti}) for $h\equiv 0$. Then we have
$$
X(t)-\xi(t)=\int_0^t\sigma(X(s))\circ dw(s)
+\int_0^t(b(X(s))-b(\xi(s)))ds-\int_0^t(dK(s)-\eta(s)ds).
$$
Set $\psi(x):=1-e^{-|x|^2}$,
 $\psi_i(x):=\frac{\p}{\p x_i}\psi(x)$,
 $\psi_{i,l}(x):=\frac{\p}{\p x_l}\psi_i(x)$.
Then $\psi\in C_b^2(\mR^m)$ and there exists $C>0$ such that
$|\psi_i(x)|=|2x_i e^{-|x|^2}|\leq C\psi(x)$ and $|\psi_{i,l}(x)|\leq C\psi(x)$.
Set again
$$G(t):=X(t)-\xi(t).$$
Since $\<X(s)-\xi(s),dK(s)-\eta(s)ds\>\geq 0$, we can have
\ce
\psi(G(t))&=&\int_0^t\psi_i(G(s))\sigma_k^i(X(s))\circ dw^k(s)\\
&&+\int_0^t\psi_i(G(s))(b^i(X(s))-b^i(\xi(s)))ds\\
&&-2\int_0^te^{-|G(s)|^2}\<X(s)-\xi(s),dK(s)-\eta(s)ds\>\\
&\leq& \int_0^t\psi_i(G(s))\sigma_k^i(X(s))\circ dw^k(s)+C\int_0^t\psi(G(s))ds\\
&:=&I+C\int_0^t\psi(G(s))ds.
\de
Let $\sigma_{k,l}^i(x):=\frac{\p}{\p x_l}\sigma_k^i(x)$,
then
\ce
I&=&\int_0^t\psi_i(G(s))\sigma_k^i(X(s))\circ dw^k(s)\\
&=&\psi_i(G(t))\sigma_k^i(X(t))w^k(t)
-\int_0^tw^k(s) \circ d[\psi_i(G(s))\sigma_k^i(X(s))]\\
&=&\psi_i(G(t))\sigma_k^i(X(t))w^k(t)
-\int_0^t[w^k(s) \psi_i(G(s))]\circ d\sigma_k^i(X(s))\\
&&-\int_0^t[w^k(s) \sigma_k^i(X(s))]\circ d(\psi_i(G(s)))\\
&=&\psi_i(G(t))\sigma_k^i(X(t))w^k(t)\\
&&-\int_0^t\psi_i(G(s)) \sigma_{k,l}^i(X(s)) \sigma_m^l(X(s))w^k(s)\circ dw^m(s)\\
&&-\int_0^t\psi_i(G(s))\sigma_{k,l}^i(X(s))b^l(X(s))w^k(s) ds\\
&&+\int_0^t\psi_i(G(s))\sigma_{k,l}^i(X(s))w^k(s)dK^l(s)\\
&&-\int_0^t\psi_{i,l}(G(s))\sigma_k^i(X(s)) \sigma_m^l(X(s))w^k(s) \circ dw^m(s)\\
&&-\int_0^t\psi_{i,l}(G(s))\sigma_k^i(X(s)) w^k(s)[b^l(X(s))-b^l(\xi(s))]ds\\
&&+\int_0^t\psi_{i,l}(G(s))\sigma_k^i(X(s)) w^k(s) [dK^l(s)-\eta(s)ds]\\
&:=&I_{1}+I_{2}+I_{3}+I_{4}+I_{5}+I_{6}+I_{7}.
\de
Obviously,
\ce
I_1\leq C\|w\|_T,~~~I_3\leq C\|w\|_T, ~~~I_4 \leq C\|w\|_T |K|_T,
\de
\ce
I_6 \leq C\int_0^t\psi(G(s))\|w\|_T ds \leq C\|w\|_T
\de
and
\ce
I_7&\leq& C\int_0^t|w(s)|d|K|_s^0+C\int_0^t|w(s)||\eta(s)|ds\\
&\leq& C\|w\|_T|K|_T+C\|w\|_T.
\de
We need to estimate $I_2$ and $I_6$. Clearly,

\ce
I_2&=&-\int_0^t\psi_i(G(s))\sigma_{k,l}^i(X(s))
\sigma_m^l(X(s))\circ d \zeta^{km}(s)\\
&=&-\int_0^t\psi_i(G(s))\sigma_{k,l}^i(X(s))\sigma_m^l(X(s)) d\zeta^{km}(s)\\
&&-\frac{1}{2}\int_0^t\frac{\partial}{\partial x_n}
[\sigma_{k,l}^i\sigma_m^l](X(s))\psi_i(G(s))
\sigma_{\a}^n(X(s))w^k(s)\d^{\a m}ds\\
&&-\frac{1}{2}\int_0^t[\sigma_{k,l}^i\sigma_m^l](X(s))\psi_{i,j}(G(s))\sigma_{\a}^j(X(s))w^k(s)\d^{\a m}ds\\
&\leq &-\int_0^t\psi_i(G(s))\sigma_{k,l}^i(X(s))\sigma_m^l(X(s)) d\zeta^{km}(s)
+C\|w\|_T.
\de
and
\ce
I_5&=&-\int_0^t\psi_{i,l}(G(s))\sigma_k^i(X(s)) \sigma_m^l(X(s))\circ d\zeta^{km}(s)\\
&=&-\int_0^t\psi_{i,l}(G(s))\sigma_k^i(X(s)) \sigma_m^l(X(s)) d\zeta^{km}(s)\\
&&-\frac{1}{2}\int_0^t\frac{\partial}{\partial x_n}
[\sigma_k^i\sigma_m^l](X(s))\psi_{i,l}(G(s))
\sigma_{\a}^n(X(s))w^k(s)\d^{\a m}ds\\
&&-\frac{1}{2}\int_0^t[\sigma_k^i\sigma_m^l](X(s))\psi_{i,l,j}(G(s))\sigma_{\a}^j(X(s))w^k(s)\d^{\a m}ds\\
&\leq &-\int_0^t\psi_{i,l}(G(s))\sigma_k^i(X(s)) \sigma_m^l(X(s)) d\zeta^{km}(s)
+C\|w\|_T.
\de

By all the above, we can get that
\ce
\psi(G(t)) &\leq& C\int_0^t\psi(G(s))ds+C\|w\|_T |K|_T+C\|w\|_T\\
&&-\int_0^t\psi_i(G(s))\sigma_{k,l}^i(X(s))\sigma_m^l(X(s)) d\zeta^{km}(s)\\
&&-\int_0^t\psi_{i,l}(G(s))\sigma_k^i(X(s)) \sigma_m^l(X(s)) d\zeta^{km}(s)\\
&:=& C\int_0^t\psi(G(s))ds+\sum_{i=1}^4A_i
\de
Obviously, for every $\e>0$, $P(\|A_i\|_T>\e|\|w\|_T<\d)\rightarrow 0 \quad as \quad\d\downarrow 0$
holds for $i=1$ and $i=2$. By Lemma \ref{lem}, for $i=3$ and $i=4$ we can have
\be
P(\|A_i\|_T>\e|\|w\|_T<\d)\rightarrow 0 \quad as \quad\d\downarrow 0. \label{am}
\ee

On the set $\{\om;\|A_i\|_T <\e,i=1,...,4\}$, we have
$$
\psi(G(t))\leq \e\exp\{CT\}\leq C \e,
$$
that is
$$
|X(t)-\xi(t)|\leq \sqrt{-\ln (1-C\e)}.
$$
Combining this with (\ref{am}), we get that
\ce
P(\|X-\xi\|_T>\e|\|w\|_T<\d)\rightarrow 0 \quad as \quad\d\downarrow 0.
\de
Finally, to see
\ce
P(|K-\eta|_T<\e|\|w-h\|_T<\d)\rightarrow 1 \quad as \quad \d\downarrow 0,
\de
it suffices to notice that
$$
K(t)-\eta(t)=X(t)-\xi(t)+\int_0^t\sigma(X(s))dw(s)+\int_0^t(b(X(s))-b(\xi(s)))ds
$$
and use Lemma \ref{fbecb2}.
\end{proof}

\end{CJK}

\end{document}